\documentclass[british,letterpaper, 10 pt, conference]{ieeeconf}
\usepackage[T1]{fontenc}
\usepackage[latin9]{inputenc}
\usepackage{mathrsfs}
\usepackage{amsmath}
\usepackage{amssymb}

\makeatletter
\IEEEoverridecommandlockouts         
\title{\LARGE \bf
Energy-Based Control of Nonlinear Infinite-Dimensional Port-Hamiltonian Systems with Dissipation
}
\author{T. Malzer, H. Rams, M. Schöberl% <-this % stops a space
\thanks{This work has been supported by the Austrian Science Fund (FWF) under
grant number P 29964-N32. All authors are with the Institute of Automatic
Control and Control Systems Technology, Johannes Kepler University Linz,
Altenbergerstrasse 66, 4040 Linz, Austria.
        {\tt\small \{tobias.malzer\_1, hubert.rams, markus.schoeberl\}@jku.at}}%
}

\usepackage{pgfplots}
 \pgfplotsset{compat=newest}
  \usetikzlibrary{plotmarks}
  \usepackage{grffile}

  \usepackage{amsthm}
  \usepackage{microtype}

  \usepackage{amsmath}
  \usepackage{mathrsfs}
  \usepackage{subcaption}
  \usepackage{graphics,xcolor}
\definecolor{P285U}{cmyk}{0.89,0.43,0.0,0.0}
\definecolor{P285U_font}{cmyk}{0.89,0.43,0.0,0.4}
\definecolor{lgray}{cmyk}{0,0,0,0.2}
\definecolor{myblue}{cmyk}{100,75,0,0}
\definecolor{jkuBlue}{RGB}{4,110,152}
\definecolor{jkuBlue}{RGB}{0,120,170}
\definecolor{jkuCyan}{RGB}{100,180,190}
\definecolor{jkuYellow}{RGB}{230,195,35}
\definecolor{jkuGrey}{RGB}{125,130,140}
\definecolor{jkuDarkGrey}{RGB}{51,51,51}
\definecolor{jkuLightGreen}{RGB}{195,215,75}
\definecolor{jkuGreen}{RGB}{115,180,85}
\definecolor{jkuPurple}{RGB}{145,75,130}
\definecolor{jkuRed}{RGB}{205,90,80}

\usepackage{mdframed}

\theoremstyle{definition}

\newtheorem{exmp}{Example}
\newtheorem{rem}{Remark}

\makeatother

\usepackage{babel}
\begin{document}
\maketitle \thispagestyle{empty} \pagestyle{empty} 
\begin{abstract}
In this paper, we consider nonlinear PDEs in a port-Hamiltonian setting
based on an underlying jet-bundle structure. We restrict ourselves
to systems with 1-dimensional spatial domain and 2nd-order Hamiltonian
including certain dissipation models that can be incorporated in the
port-Hamiltonian framework by means of appropriate differential operators.
For this system class, energy-based control by means of Casimir functionals
as well as energy balancing is analysed and demonstrated using a nonlinear
Euler-Bernoulli beam.
\end{abstract}

\section{Introduction}

The port-Hamiltonian (pH) system formulation in combination with energy-based
control schemes has turned out to be an effective tool for the description
and control of nonlinear finite-dimensional systems, see \cite{Schaft2000,Ortega2001}
for instance. Especially, for the stabilisation of, e.g., non-energy
minimal equilibria the control by interconnection based on Casimir
functionals as well as the energy balancing methodology are well established.

In the infinite-dimensional scenario, the pH system representation
is \textendash{} in contrast to the finite-dimensional scenario \textendash{}
not unique. With regard to control engineering purposes the Stokes-Dirac
approach, see \cite{Schaft2002,Gorrec2005,Jacob2012}, as well as
an approach based on jet-bundle structures, see \cite{Schoeberl2014a,Schoeberl2015e},
have turned out to be adequate frameworks. These approaches mainly
differ in the choice of variables (energy variables versus derivative
variables in the Hamiltonian). This has the consequence that for linear
partial differential equations (PDEs) the Stokes-Dirac framework is
closely connected to well-known functional-analytic methods, in particular
to the theory of strongly continuous semigroups \cite{Curtain1995},
whereas the jet-bundle approach is very well suited for systems that
allow for a variational characterisation including nonlinear systems.
One of the major benefits of infinite-dimensional pH system formulations
is that they provide a consistent framework for the development of
(finite-dimensional) boundary controllers which are of great practical
relevance. In view of this, the two most important schemes that have
already been extended to the infinite-dimensional setting are the
energy-Casimir method (see, e.g., \cite{Macchelli2004} for the Stokes-Dirac
approach and \cite{Siuka2011a,Rams2017a} for the jet-bundle framework)
as well as the control by energy balancing (see \cite{Macchelli2017}).
Note that in \cite{Macchelli2004,Macchelli2017} as well as in \cite{Siuka2011a,Rams2017a}
only linear mechanical structures have been considered to show the
applicability of the proposed control schemes. It is worth stressing
that, similar to the finite-dimensional case, both methods can be
used to stabilise, e.g., non-energy minimal rest positions; however,
the energy-Casimir method yields dynamic controllers whereas energy
balancing usually leads to static control laws. 

In this paper, we focus on the description of (nonlinear) mechanical
systems on 1-dimensional spatial domains formulated within the jet-bundle
approach. To demonstrate the capability of this approach, we refer
to \cite{Siuka2011a,Rams2017a} for systems with 1st- and 2nd-order
Hamiltonian densities. The order of the Hamiltonian density is basically
responsible for the number of boundary-port categories that can be
introduced. Therefore, in what follows, we confine ourselves to 2nd-order
Hamiltonian densities. In general, the considered jet-bundle approach,
where the Hamiltonian depends on derivative variables, can be divided
in the non-differential case and the more general differential-operator
case with respect to the interconnection and damping maps, see \cite{Schoeberl2014a}.
In this setup, for mechanical systems the interconnection map can
usually be chosen as the canonical map and nontrivial dissipation
maps can be used to include certain damping models. It should be stressed
that within the jet-bundle approach boundary ports are a consequence
of the derivative variables in the Hamiltonian; however, using differential
operators for the dissipation mapping leads to modified boundary-port
relations. Consequently, the introduction of boundary ports is not
straightforward and depends mainly on the concrete operators under
consideration.

Therefore, the main contributions of this paper are as follows: i)
we study the impact of certain differential operators \textendash{}
used as dissipation mappings \textendash{} on the boundary ports of
a nonlinear infinite-dimensional pH-system in the jet-bundle approach,
see Section \ref{sec:Infinite-Dimensional-PH-Systems}; ii) we propose
a control scheme based on the energy-Casimir method that exploits
a certain boundary-output assignment in the pH-framework including
differential operators (e.g. nontrivial dissipation maps) and 2nd-order
Hamiltonian densities, see Section \ref{sec:Casimir-Control}; iii)
the control by energy balancing (EBC) is investigated for nonlinear
PDEs within the jet-bundle framework mainly exploiting the geometric
properties of the corresponding boundary operators, see Section \ref{sec:Energy-Balancing-Control}.
To show the applicability of the proposed theory including the control
strategies, the example of a nonlinear Euler-Bernoulli beam subject
to structural damping is considered.

\section{Notation and Preliminaries\label{sec:Notation-and-Preliminaries}}

In this paper, we make heavy use of differential-geometric methods,
with a notation similar to \cite{Giachetta1997}. To keep the formulas
short and readable, we apply tensor notation and Einsteins convention
on sums. The symbols $\wedge$, $\rfloor$ and $\mathrm{d}$ denote
the exterior (wedge) product, the natural contraction between tensor
fields and the exterior derivative, respectively. It should be noted
that the use of pull-back bundles is omitted for ease of presentation.
Furthermore, the ranges of the used indices are not indicated when
they are clear from the context. The set of all smooth functions on
an arbitrary manifold $\mathcal{M}$ is denoted by $C^{\infty}(\mathcal{M})$.

In the following, we investigate PDEs with 1-dimensional spatial domain.
To be able to distinguish between dependent and independent coordinates,
we introduce bundle structures. Let us consider the bundle $\pi:\mathcal{E}\rightarrow\mathcal{B}$,
where $\pi$ is a surjective submersion \textendash{} called projection
\textendash{} from the total manifold $\mathcal{E}$ to the base manifold
$\mathcal{B}$. Since we confine ourselves to 1-dimensional spatial
domains, $\mathcal{B}$ only possesses the independent (spatial) coordinate
$z^{1}$. Note that $\partial\mathcal{B}$ represents the boundary
of the manifold $\mathcal{B}$, and the restriction of a mathematical
expression to $\partial\mathcal{B}$ is indicated with $(\cdot)|_{\partial\mathcal{B}}$.
Furthermore, the total manifold $\mathcal{E}$ comprises the coordinates
$(z^{1},x^{\alpha})$ with $\alpha=1,\ldots,n$, where $x^{\alpha}$
denote the dependent coordinates. Next, let us introduce derivative
coordinates (jet variables). To this end, we consider the (higher-order)
jet manifold $\mathcal{J}^{r}(\mathcal{E})$, possessing the coordinates
$(z^{1},x^{\alpha},x_{1}^{\alpha},\ldots,x_{1\ldots1}^{\alpha})$,
where for instance $x_{111}^{\alpha}$ denotes the 3rd-order derivative
coordinate, i.e. the 3rd derivative of $x^{\alpha}$ with respect
to the independent coordinate $z^{1}$. Here, $r$ denotes the highest
occuring order of derivatives and, exemplarily, for the 4th jet manifold
$\mathcal{J}^{4}(\mathcal{E})$, we have the coordinates $(z^{1},x^{\alpha},x_{1}^{\alpha},x_{11}^{\alpha},x_{111}^{\alpha},x_{1111}^{\alpha})$.

Furthermore, we introduce the tangent bundle $\tau_{\mathcal{E}}:\mathcal{T}(\mathcal{E})\rightarrow\mathcal{E}$
possessing the coordinates $(z^{1},x^{\alpha},\dot{z}^{1},\dot{x}^{\alpha})$,
together with the abbreviations $\partial_{1}=\partial/\partial z^{1}$
and $\partial_{\alpha}=\partial/\partial x^{\alpha}$ denoting the
fibre bases of the bundle. Of special interest is the so-called vertical
tangent bundle $\nu:\mathcal{V}(\mathcal{E})\rightarrow\mathcal{E}$,
which is a subbundle of $\tau_{\mathcal{E}}$ and equipped with the
coordinates $(z^{1},x^{\alpha},\dot{x}^{\alpha})$. By means of the
total derivative $d_{1}=\partial_{1}+x_{1}^{\alpha}\partial_{\alpha}+x_{11}^{\alpha}\partial_{\alpha}^{1}+x_{111}^{\alpha}\partial_{\alpha}^{11}+x_{1111}^{\alpha}\partial_{\alpha}^{111}+\ldots$,
we are able to introduce the $n$th prolongation of a vertical vector
field $v=v^{\alpha}\partial_{\alpha}$, where for instance the $2$nd
prolongation reads as $j^{2}\left(v\right)=v^{\alpha}\partial_{\alpha}+d_{1}(v^{\alpha})\partial_{\alpha}^{1}+d_{11}(v^{\alpha})\partial_{\alpha}^{11}$,
with the abbreviations $\partial_{\alpha}^{1}=\partial/\partial x_{1}^{\alpha}$,
$\partial_{\alpha}^{1\ldots1}=\partial/\partial x_{1\ldots1}^{\alpha}$
and $d_{11}=d_{1}\circ d_{1}$ denoting the repeated total derivative.

In what follows, we also need some futher important structures in
order to be able to consider one-forms. The cotangent bundle $\tau_{\mathcal{E}}^{*}=\mathcal{T}^{*}\left(\mathcal{E}\right)\rightarrow\mathcal{E}$,
where we have the coordinates $(z^{1},x^{\alpha},\dot{z}_{1},\dot{x}_{\alpha})$
and the holonomic bases $\mathrm{d}z^{1},\mathrm{d}x^{\alpha}$, allows
us to locally define a section $\omega=\omega_{1}\mathrm{d}z^{1}+\omega_{\alpha}\mathrm{d}x^{\alpha}$,
with $\omega_{1},\omega_{\alpha}\in C^{\infty}(\mathcal{E})$, on
it. In this paper, we focus on (Hamiltonian) densities $\mathfrak{H}=\mathcal{H}\Omega$
and (Hamiltonian) functionals $\mathscr{H}=\int_{\mathcal{B}}\mathcal{H}\Omega$,
with $\mathcal{H}\in C^{\infty}(\mathcal{J}^{2}(\mathcal{E}))$ \textendash{}
i.e. on densities that may depend on derivative coordinates \textendash ,
where $\Omega=\mathrm{d}z^{1}$ denotes the volume element on $\mathcal{B}$
and $\Omega_{1}=\partial_{1}\rfloor\mathrm{d}z^{1}$ the boundary-volume
form. The bundle structure $\pi:\mathcal{E}\rightarrow\mathcal{B}$
allows us to construct some further geometric objects like the tensor
bundle $\mathcal{W}_{1}^{r}(\mathcal{E})=\mathcal{T}^{*}(\mathcal{E})\wedge\mathcal{T}^{*}(\mathcal{B})$
with a typical element $\omega_{\alpha}\mathrm{d}x^{\alpha}\wedge\mathrm{d}z^{1}$
for $\mathcal{W}_{1}^{r}(\mathcal{E})$, where $\omega_{\alpha}\in C^{\infty}(\mathcal{J}^{r}(\mathcal{E}))$
is met. Moreover, we consider $k$th-order linear differential operators
$\mathfrak{D}:\mathcal{W}_{1}^{r}(\mathcal{E})\rightarrow\mathcal{V}(\mathcal{E})$
serving as a map of an element $\mathcal{W}_{1}^{r}(\mathcal{E})$
of jet order $r$ to an element $\mathcal{V}(\mathcal{E})$ of jet
order $r+k$.

Further, we introduce the so-called horizontal exterior derivative
$\mathrm{d}_{h}$, meeting $\mathrm{d}_{h}(\omega)=\mathrm{d}z^{1}\wedge d_{1}(\omega)$
for a form $\omega:\mathcal{J}^{r}(\mathcal{E})\rightarrow\mathcal{T}^{*}\left(\mathcal{J}^{r}(\mathcal{E})\right)$,
to be able to make use of Stokes theorem, see \cite{Giachetta1997}
and \cite[Appendix]{Schoeberl2014a} for more details.

\section{Infinite-Dimensional PH-Systems\label{sec:Infinite-Dimensional-PH-Systems}}

In this section, we extend the pH-framework presented in \cite[Definition 4]{Schoeberl2014a}
to 2nd-order Hamiltonian densities with the restriction on 1-dimensional
spatial domains by following the findings of \cite{Rams2017a}. This
framework is mainly based on the underlying jet-bundle structures
of the PDEs under consideration and makes heavy use of a certain power-balance
relation. Moreover, we solely investigate systems with boundary in-
and outputs, i.e. in-domain in- and outputs are not addressed within
this paper.

Let $\mathfrak{H}$ be a 2nd-order Hamiltonian, i.e. $\mathcal{H}\in C^{\infty}(\mathcal{J}^{2}(\mathcal{E}))$,
then a pH system formulation using differential operators is
\begin{equation}
\dot{x}=(\mathfrak{J}-\mathfrak{R})(\delta\mathfrak{H}),\label{eq:iPCHD_Operator}
\end{equation}
including appropriate boundary conditions. Here, $\mathfrak{J},\,\mathfrak{R}$
are $r$th-order linear vector-valued differential operators and serve
as maps $\mathfrak{J},\mathfrak{R}:\mathcal{W}_{1}^{4}(\mathcal{E})\rightarrow\mathcal{V}(\mathcal{E})$
with the following properties. The interconnection operator $\mathfrak{J}$
characterises the internal power flow, and enjoys the property of
formal skew-adjointness, i.e. $\mathfrak{J}(\eta)\rfloor\varpi=-\mathfrak{J}(\varpi)\rfloor\eta+\mathrm{d}_{h}(\mathfrak{j})$,
with $\mathfrak{j}=\mathfrak{j}^{1}\Omega_{1}$ and $\eta,\varpi\in\mathcal{W}_{1}^{4}(\mathcal{E})$.
In addition, $\mathfrak{R}$ comprises dissipative effects and is
described by a formally self-adjoint, non-negative operator meeting
$\mathfrak{R}(\eta)\rfloor\varpi=\mathfrak{R}(\varpi)\rfloor\eta+\mathrm{d}_{h}(\mathfrak{r})$
with $\mathfrak{r}=\mathfrak{r}^{1}\Omega_{1}$, and the non-negativity
relation $\mathfrak{R}(\eta)\rfloor\eta\geq0$. Furthermore, for 2nd-order
Hamiltonian densities, the variational derivative corresponds to $\delta\mathfrak{H}=\delta_{\alpha}\mathcal{H}\mathrm{d}x^{\alpha}\wedge\Omega$
with $\delta_{\alpha}(\cdot)=\partial_{\alpha}(\cdot)-d_{1}(\partial_{\alpha}^{1}(\cdot))+d_{11}(\partial_{\alpha}^{11}(\cdot))$.
It is of particular interest how the Hamiltonian functional $\mathscr{H}$
evolves along solutions of the system (\ref{eq:iPCHD_Operator}) (well-posedness
provided). For (\ref{eq:iPCHD_Operator}) the formal change can be
deduced to the balance relation
\begin{equation}
\dot{\mathscr{H}}=\int_{\mathcal{B}}(\mathfrak{J}-\mathfrak{R})(\delta\mathfrak{H})\rfloor\delta\mathfrak{H}+(\dot{x}\rfloor\delta^{\partial,1}\mathfrak{H}+\dot{x}_{1}\rfloor\delta^{\partial,2}\mathfrak{H})|_{\partial\mathcal{B}},\label{eq:h_p_operator_case}
\end{equation}
which states a power-balance relation if $\mathscr{H}$ represents
the total energy of the system. From a control engineering point of
view, it is of special interest to introduce power ports by means
of (\ref{eq:h_p_operator_case}), where we basically exploit both
the boundary operators $\delta^{\partial,1}\mathfrak{H}=(\partial_{\alpha}^{1}\mathcal{H}-d_{1}(\partial_{\alpha}^{11}\mathcal{H}))\mathrm{d}x^{\alpha}\wedge\Omega_{1}$
and $\delta^{\partial,2}\mathfrak{H}=\partial_{\alpha}^{11}\mathcal{H}\mathrm{d}x_{1}^{\alpha}\wedge\Omega_{1}$.\begin{rem}\label{rem:boundary_ports_diff_op_case}It
is worth stressing that boundary ports can be generated in two different
ways. On the one hand, they can be a direct consequence of the jet
variables that may occur in the Hamiltonian, cf. $\delta^{\partial,1}$
and $\delta^{\partial,2}$; on the other hand, boundary terms can
also stem from the differential operators $\mathfrak{J}$ and $\mathfrak{R}$
\textendash{} depending on their structure \textendash{} which influence
the boundary ports that are due to the jet variables, or even create
additional ports. \end{rem}\begin{rem} In the Stokes-Dirac scenario
no jet variables occur in $\mathcal{H}$ because energy variables
are used; therefore, boundary ports solely stem from the interconnection
operator $\mathfrak{J}$, see \cite{Schaft2002} for instance.\end{rem}Due
to the fact that the general introduction of boundary ports is not
possible in the operator case, we first consider the so-called non-differential
operator case, where the introduction of boundary ports is straightforward.
Second, to highlight the fact that $\mathfrak{J}$ and $\mathfrak{R}$
modify the boundary terms, we additionally investigate the impact
of specific dissipation operators in the course of Subsection \ref{subsec:Specific-Operators-for-R}.

\subsection{Non-Differential Operator Case}

In the non-differential operator case the differential operators $\mathfrak{J}$
and $\mathfrak{R}$ degenerate to bounded linear mappings $\mathcal{J},\mathcal{R}$.
Hence, in this case (\ref{eq:iPCHD_Operator}) reads in local coordinates
as
\begin{equation}
\dot{x}^{\alpha}=(\mathcal{J}^{\alpha\beta}-\mathcal{R}^{\alpha\beta})\delta_{\beta}\mathcal{H},\quad\alpha,\beta=1,\ldots,n,\label{eq:pH_non_differential_operator}
\end{equation}
where the skew-symmetric interconnection map $\mathcal{J}$ meets
$\mathcal{J}^{\alpha\beta}=-\mathcal{J}^{\beta\alpha}\in C^{\infty}(\mathcal{J}^{4}(\mathcal{E}))$,
and the positive semi-definiteness and symmetry of $\mathcal{R}$
implies $\mathcal{R}^{\alpha\beta}=\mathcal{R}^{\beta\alpha}\in C^{\infty}(\mathcal{J}^{4}(\mathcal{E}))$
and $\left[\mathcal{R}^{\alpha\beta}\right]\geq0$ for the coefficient
matrix. Consequently, the power-balance relation (\ref{eq:h_p_operator_case})
reduces to
\begin{multline}
\dot{\mathscr{H}}=-\int_{\mathcal{B}}\delta_{\alpha}(\mathcal{H})\mathcal{R}^{\alpha\beta}\delta_{\beta}(\mathcal{H})\mathrm{d}z^{1}+\ldots\\
+(\dot{x}^{\alpha}\delta_{\alpha}^{\partial,1}\mathcal{H}+\dot{x}_{1}^{\alpha}\delta_{\alpha}^{\partial,2}\mathcal{H})|_{\partial\mathcal{B}}.\label{eq:h_p_non_diff_op}
\end{multline}
Here, the first term denotes the distributed dissipation on the domain,
while the second part allows us to define power ports with
\begin{eqnarray}
\delta_{\alpha}^{\partial,1}\mathcal{H} & = & \partial_{\alpha}^{1}\mathcal{H}-d_{1}(\partial_{\alpha}^{11}\mathcal{H})\label{eq:boundary_operators}\\
\delta_{\alpha}^{\partial,2}\mathcal{H} & = & \partial_{\alpha}^{11}\mathcal{H}.\nonumber 
\end{eqnarray}
\begin{rem}From (\ref{eq:boundary_operators}) it is obvious that
in the non-differential operator case the boundary ports solely stem
from the derivative variables occuring in the Hamiltonian density.\end{rem}Henceforth,
we suppose that the boundary $\partial\mathcal{B}=\{0,L\}$ can be
divided into an unactuated part $\partial\mathcal{B}_{u}=0$ and a
fully actuated part $\partial\mathcal{B}_{a}=L$. For $\partial\mathcal{B}_{u}$,
we assume that $\dot{x}^{\alpha}\delta_{\alpha}^{\partial,1}\mathcal{H}|_{\partial\mathcal{B}_{u}}=0$
as well as $\dot{x}_{1}^{\alpha}\delta_{\alpha}^{\partial,2}\mathcal{H}|_{\partial\mathcal{B}_{u}}=0$,
i.e., no power exchange takes place through the unactuated boundary
part. For $\partial\mathcal{B}_{a}$, we set $\dot{x}^{\alpha}\delta_{\alpha}^{\partial,1}\mathcal{H}|_{\partial\mathcal{B}_{a}}=\hat{u}^{\hat{\xi}}\hat{y}_{\hat{\xi}}$
as well as $\dot{x}_{1}^{\alpha}\delta_{\alpha}^{\partial,2}\mathcal{H}|_{\partial\mathcal{B}_{a}}=\check{u}^{\check{\xi}}\check{y}_{\check{\xi}}$
with the collocated boundary pairs $(\hat{u}^{\hat{\xi}},\hat{y}_{\hat{\xi}})$
and $(\check{u}^{\check{\xi}},\check{y}_{\check{\xi}})$ including
the index ranges $\hat{\xi}=1,\ldots,\hat{m}$ and $\check{\xi}=1,\ldots,\check{m}.$
By setting
\begin{equation}
\begin{array}{ccc}
\hat{B}_{\alpha\hat{\xi}}\hat{u}^{\hat{\xi}} & = & \delta_{\alpha}^{\partial,1}\mathcal{H}|_{\partial\mathcal{B}_{a}},\\
\hat{y}_{\hat{\xi}} & = & \hat{B}_{\alpha\hat{\xi}}\dot{x}^{\alpha}|_{\partial\mathcal{B}_{a}},
\end{array}\;\begin{array}{ccc}
\check{B}_{\alpha\check{\xi}}\check{u}^{\check{\xi}} & = & \delta_{\alpha}^{\partial,2}\mathcal{H}|_{\partial\mathcal{B}_{a}},\\
\check{y}_{\check{\xi}} & = & \check{B}_{\alpha\check{\xi}}\dot{x}_{1}^{\alpha}|_{\partial\mathcal{B}_{a}},
\end{array}\label{eq:output_param}
\end{equation}
we assign the roles of the inputs and outputs in (\ref{eq:h_p_non_diff_op}),
which is of course not unique, see \cite{Schoeberl2008a}. To highlight
the benefits of the pH system representation for the non-differential
operator scenario, in the appendix we study the Examples \ref{ex:lin_beam}
and \ref{ex:nonlin_beam} where the transversal and the longitudinal
deflection of a linear and a nonlinear beam structure are investigated.

\subsection{Specific Operators for $\mathfrak{R}$\label{subsec:Specific-Operators-for-R}}

The objective of this section is to study the transversal as well
as the longitudinal deflection of a nonlinear beam model subject to
structural damping. To this end, we introduce two linear differential
operators and investigate their impact on the power-balance equation,
more specifically how they affect the geometric boundary ports.

First, we consider the 2nd-order formally self-adjoint operator $\mathfrak{R}_{A}$,
locally expressed by
\[
\mathfrak{R}_{A}(\eta)=d_{1}(\mathfrak{R}_{A}^{\alpha\beta}d_{1}(\eta_{\beta}))\partial_{\alpha},
\]
with $\mathfrak{R}_{A}^{\alpha\beta}=\mathfrak{R}_{A}^{\beta\alpha}\in C^{\infty}(\mathcal{B})$,
i.e. the coefficients may in general depend on the spatial coordinate.
Furthermore, we are able to deduce the important relation
\begin{equation}
\mathfrak{R}_{A}(\eta)\rfloor\eta=\mathrm{d}_{h}(\mathfrak{R}_{A}^{\alpha\beta}d_{1}(\eta_{\beta})\eta_{\alpha}\Omega_{1})-d_{1}(\eta_{\beta})\mathfrak{R}_{A}^{\alpha\beta}d_{1}(\eta_{\alpha})\Omega,\label{eq:second_order_op_decomp}
\end{equation}
indicating that $\mathfrak{R}_{A}$ influences the (geometric) boundary
ports (first term in (\ref{eq:second_order_op_decomp})) as well as
the domain conditions (second term in (\ref{eq:second_order_op_decomp})).
The non-negativity of $\mathfrak{R}_{A}$ follows if $d_{1}(\eta_{\beta})\mathfrak{R}_{A}^{\alpha\beta}d_{1}(\eta_{\alpha})\Omega\leq0$
is met, i.e. the coefficient matrix $[\mathfrak{R}_{A}^{\alpha\beta}]$
must be symmetric as well as negative semi-definite.

Next, to be able to describe the structural damping for the vertical
deflection, we exploit the 4th-order operator
\[
\mathfrak{R}_{B}(\eta)=d_{11}(\mathfrak{R}_{B}^{\alpha\beta}d_{11}(\eta_{\beta}))\partial_{\alpha},
\]
with $\mathfrak{R}_{A}^{\alpha\beta}\in C^{\infty}\left(\mathcal{B}\right)$,
where $\mathfrak{R}_{A}^{\alpha\beta}=\mathfrak{R}_{A}^{\beta\alpha}$
must be met for the self-adjointness of the operator. Similar to (\ref{eq:second_order_op_decomp}),
we rewrite the expression $\mathfrak{R}_{B}\left(\eta\right)\rfloor\eta$
according to
\begin{multline}
\mathfrak{R}_{B}(\eta)\rfloor\eta=d_{11}(\eta_{\alpha})\mathfrak{R}_{B}^{\alpha\beta}d_{11}(\eta_{\beta})\Omega+\ldots\\
+\mathrm{d}_{h}((d_{1}(\mathfrak{R}_{B}^{\alpha\beta}d_{11}(\eta_{\alpha}))\eta_{\beta}-\mathfrak{R}_{B}^{\alpha\beta}d_{11}(\eta_{\alpha})d_{1}(\eta_{\beta}))\Omega_{1}),\label{eq:fourth_order_op_decomp}
\end{multline}
and find that $\mathfrak{R}_{B}$ also affects the boundary ports
as well as the domain conditions. Note that the non-negativity of
$\mathfrak{R}_{B}$ follows if $d_{11}(\eta_{\alpha})\mathfrak{R}_{B}^{\alpha\beta}d_{11}(\eta_{\beta})\geq0$
is satisfied, i.e., the coefficient matrix $[\mathfrak{R}_{B}^{\alpha\beta}]$
has to be positive semi-definite to meet this requirement.

Now, having the preceding findings at hand, we are able to investigate
the following example.

\begin{exmp}[Nonlinear 2nd-order beam structure with structural damping]\label{ex:nonlin_beam_diss}We
consider the following set of nonlinear PDEs, where \begin{subequations}\label{eq:nonlin_beam_diss_eom_0}
\begin{multline}
\rho A\ddot{w}^{1}=\frac{3}{2}EAw_{11}^{1}\left(w_{1}^{1}\right)^{2}+EAw_{11}^{2}w_{1}^{1}+\ldots\\
+EAw_{1}^{2}w_{11}^{1}-EIw_{1111}^{1}-\alpha_{1}\dot{w}_{1111}^{1}\label{eq:nonlin_beam_diss_eom_1}
\end{multline}
describes the transversal deflection $w^{1}$ and
\begin{equation}
\rho A\ddot{w}^{2}=EAw_{11}^{2}+EAw_{1}^{1}w_{11}^{1}+\alpha_{2}\dot{w}_{11}^{2}\label{eq:nonlin_beam_diss_eom_2}
\end{equation}
\end{subequations}the longitudinal deflection $w^{2}$ of a nonlinear
beam, with $E,I,\rho,A>0$ as material parameters and $\alpha_{1},\alpha_{2}>0$
as damping coefficients; see \cite{Liu1998} for instance for a similar
model. Note that (\ref{eq:nonlin_beam_diss_eom_0}) follows from the
calculus of variations \cite{Meirovitch1997} where nonlinear strain-deflection
relations in the sense of von Kármán \cite{Kugi2001} and the Euler-Bernoulli
hypothesis have been used; eventually, (linear) structural damping
has been included via the damping terms $-\alpha_{1}\dot{w}_{1111}^{1}$
and $\alpha_{2}\dot{w}_{11}^{2}$. If we introduce the generalised
momenta $p_{1}=\rho A\dot{w}^{1}$, $p_{2}=\rho A\dot{w}^{2}$ and
set the Hamiltonian density to
\begin{multline*}
\mathcal{H}=\frac{1}{2\rho A}\left((p_{1})^{2}+(p_{2})^{2}\right)+\ldots\\
+\frac{1}{2}EA((w_{1}^{2})^{2}+\frac{1}{4}(w_{1}^{1})^{4}+w_{1}^{2}(w_{1}^{1})^{2})+\frac{1}{2}EI(w_{11}^{1})^{2},
\end{multline*}
we find that the interconnection tensor
\[
\mathcal{J}=\left[\begin{array}{cccc}
0 & 0 & 1 & 0\\
0 & 0 & 0 & 1\\
-1 & 0 & 0 & 0\\
0 & -1 & 0 & 0
\end{array}\right]
\]
and the linear dissipation operator
\[
\mathfrak{R}=\left[\begin{array}{cccc}
0 & 0 & 0 & 0\\
0 & 0 & 0 & 0\\
0 & 0 & d_{11}(\alpha_{1}d_{11}(\cdot)) & 0\\
0 & 0 & 0 & -d_{1}(\alpha_{2}d_{1}(\cdot))
\end{array}\right]
\]
yield an appropriate pH system representation for (\ref{eq:nonlin_beam_diss_eom_0})
according to $\dot{x}=(\mathcal{J}-\mathfrak{R})(\delta\mathfrak{H})$.
Note that $\mathfrak{R}$ comprises the operators $\mathfrak{R}_{A}$
and $\mathfrak{R}_{B}$, i.e. $\mathfrak{R}=\mathfrak{R}_{A}+\mathfrak{R}_{B}$
with $\mathfrak{R}_{A}^{44}=-\alpha_{2}$, $\mathfrak{R}_{B}^{33}=\alpha_{1}$,
$\mathfrak{R}_{A}^{ij}=0$ for $i,j\neq4$ and $\mathfrak{R}_{B}^{kl}=0$
for $k,l\neq3$. It is worth stressing that the evaluation of the
power-balance relation (\ref{eq:h_p_operator_case}) is not as straightforward
as in the non-differential operator case. Since this fact is a key
point of this paper, we explain the following steps in detail. To
this end, we consider the expression
\[
-\mathfrak{R}\left(\delta\mathfrak{H}\right)=\left[\begin{array}{c}
0\\
0\\
-d_{11}(\alpha_{1}d_{11}(\frac{p_{1}}{\rho A}))\\
d_{1}(\alpha_{2}d_{1}(\frac{p_{2}}{\rho A}))
\end{array}\right],
\]
which plays an important role in evaluating (\ref{eq:h_p_operator_case});
therefore, the power-balance relation reads as
\begin{multline}
\dot{\mathscr{H}}=\int_{\mathcal{B}}(-d_{11}(\alpha_{1}\dot{w}_{11}^{1})\dot{w}^{1}+d_{1}(\alpha_{2}\dot{w}_{1}^{2})\dot{w}^{2})\mathrm{d}z^{1}+\ldots\\
+(\dot{w}^{1}\tilde{Q}+\dot{w}^{2}\tilde{N}+\dot{w}_{1}^{1}\tilde{M})|_{\partial\mathcal{B}_{a}},\label{eq:h_p_nonlin_beam_diss}
\end{multline}
where the shear force, the normal force and the bending moment follow
to
\begin{align}
\tilde{Q} & =\frac{1}{2}EA(w_{1}^{1})^{3}+EAw_{1}^{2}w_{1}^{1}-EIw_{111}^{1},\nonumber \\
\tilde{N} & =EAw_{1}^{2}+\frac{1}{2}EA(w_{1}^{1})^{2},\label{eq:internal_force_boundary_operators}\\
\tilde{M} & =EIw_{11}^{1},\nonumber 
\end{align}
 by applying (\ref{eq:boundary_operators}). However, if we use (\ref{eq:second_order_op_decomp})
and (\ref{eq:fourth_order_op_decomp}), the power-balance relation
(\ref{eq:h_p_nonlin_beam_diss}) can be rewritten as
\begin{multline}
\dot{\mathscr{H}}=-\int_{\mathcal{B}}(\alpha_{1}(\dot{w}_{11}^{1})^{2}+\alpha_{2}(\dot{w}_{1}^{2})^{2})\mathrm{d}z^{1}+\ldots\\
+(\dot{w}^{1}\breve{Q}+\dot{w}^{2}\breve{N}+\dot{w}_{1}^{1}\breve{M})|_{\partial\mathcal{B}_{a}},\label{eq:h_p_nonlin_beam_diss_1}
\end{multline}
where it becomes clear that $\alpha_{1},\alpha_{2}>0$ ensures the
non-negativity of $\mathfrak{R}$ which \textendash{} provided that
no power flow takes place via $\partial\mathcal{B}_{a}$ \textendash{}
guarantees the energy dissipation according to the structural damping
models under consideration. Furthermore, comparing the expressions
\[
\begin{array}{ccc}
\breve{Q}=\tilde{Q}-\alpha_{1}\dot{w}_{111}^{1}, & \breve{N}=\tilde{N}+\alpha_{2}\dot{w}_{1}^{2}, & \breve{M}=\tilde{M}+\alpha_{1}\dot{w}_{11}^{1}\end{array}
\]
with (\ref{eq:internal_force_boundary_operators}) highlights the
impact of the dissipation operator $\mathfrak{R}$ on the (geometric)
boundary-port relations. In accordance with Remark \ref{rem:boundary_ports_diff_op_case},
we find that $\breve{Q}$, $\breve{N}$ and $\breve{M}$ stem, on
the one hand, from the application of the boundary operators (\ref{eq:boundary_operators})
which are a consequence of the jet variables in $\mathcal{H}$ and,
on the other hand, from the impact of the differential operators.
This is important to realise since we use the shear force, the normal
force and the bending moment at $z^{1}=L$ as manipulated variables
for the controller designs in the following two sections. Therefore,
the collocated inputs and outputs are $\hat{u}^{1}=\breve{Q}$, $\hat{u}^{2}=\breve{N}$,
$\check{u}^{1}=\breve{M}$ and $\hat{y}_{1}=\dot{w}^{1}$, $\hat{y}_{2}=\dot{w}^{2}$,
$\check{y}_{1}=\dot{w}_{1}^{1}$, respectively. Consequently, the
boundary maps result in $\hat{B}_{\alpha\hat{\xi}}=\delta_{\alpha\hat{\xi}}$
and $\check{B}_{\alpha\check{\xi}}=\delta_{\alpha\check{\xi}}$, together
with the Kronecker-Delta symbol meeting $\delta_{\alpha\beta}=1$
for $\alpha=\beta$ and $\delta_{\alpha\beta}=0$ for $\alpha\neq\beta$.\end{exmp}

Ex. \ref{ex:nonlin_beam_diss} clearly highlights that differential
operators have a strong impact on the power-balance relation and therefore
on the introduction of the associated power ports.

The aim of the next two sections is to stabilise a certain rest position
of the nonlinear beam structure including damping, see Ex. \ref{ex:nonlin_beam_diss}.
To this end, we shall use a dynamic pH-controller exploiting Casimir
functionals as well as a static control law designed by means of energy
balancing.

\section{Energy-Based Control by means of Casimir Functionals\label{sec:Casimir-Control}}

In this section, we extend the results and findings for the energy-Casimir
method applied to systems with 2nd-order Hamiltonians in the non-differential
operator scenario, see \cite{Rams2017a}, to the differential operator
setting with a certain input and output assignment for the boundary
ports. For detailed informations concerning the 1st-order and 2nd-order
case, we refer to \cite{Schoeberl2013a} and \cite{Rams2017a}, respectively.
Since the introduction of boundary ports in the differential-operator
scenario is not straightforward, cf. Ex. \ref{ex:nonlin_beam_diss},
we suppose that, independently of the concrete differential operators,
the outputs of the considered pH-system can be parameterised via (\ref{eq:output_param}).
Note that the parameterisation of the corresponding collocated inputs
strongly depends on the involved differential operators and cannot
be stated in a general form. However, this assumption enables us to
extend the energy-Casimir method to the differential-operator scenario
as well.

\subsection{Interconnection (Infinite-Finite)}

In the following, we are interested in a power-conserving interconnection
of the infinite-dimensional plant (\ref{eq:iPCHD_Operator}) and a
finite-dimensional controller at the actuated boundary $\partial\mathcal{B}_{a}$
according to
\begin{equation}
\hat{u}\rfloor\hat{y}+\check{u}\rfloor\check{y}+\hat{u}_{c}\rfloor\hat{y}_{c}+\check{u}_{c}\rfloor\check{y}_{c}=0.\label{eq:power_conserving_interconnection}
\end{equation}
Here, $\hat{u}_{c}\rfloor\hat{y}_{c}$ and $\check{u}_{c}\rfloor\check{y}_{c}$
represent the collocated in- and output pairings of the controller,
which have been divided into two parts to take account of both boundary-port
categories of the plant. The structure of the dynamic pH-controller
can locally be given as
\begin{align}
\dot{x}_{c}^{\alpha_{c}} & =(J_{c}^{\alpha_{c}\beta_{c}}-R_{c}^{\alpha_{c}\beta_{c}})\partial_{\beta_{c}}H_{c}+\hat{G}_{c,\hat{i}}^{\alpha_{c}}\hat{u}_{c}^{\hat{i}}+\check{G}_{c,\check{j}}^{\alpha_{c}}\check{u}_{c}^{\check{j}}\nonumber \\
\hat{y}_{c,\hat{i}} & =\hat{G}_{c,\hat{i}}^{\alpha_{c}}\partial_{\alpha_{c}}H_{c}\quad\text{\ensuremath{\mathrm{and}}}\quad\check{y}_{c,\check{j}}=\check{G}_{c,\check{j}}^{\alpha_{c}}\partial_{\alpha_{c}}H_{c}\label{eq:pH_controller}
\end{align}
with $\alpha_{c},\beta_{c}=1,\ldots,n_{c}$, $\hat{i}=1,\ldots,\hat{m}$
and $\check{j}=1,\ldots,\check{m}$. Note that $\dot{H}_{c}$ can
be deduced to $-\partial_{\alpha_{c}}(H_{c})R_{c}^{\alpha_{c}\beta_{c}}\partial_{\beta_{c}}(H_{c})+\hat{u}_{c}^{\hat{i}}\hat{y}_{c,\hat{i}}+\check{u}_{c}^{\check{i}}\check{y}_{c,\check{i}}$;
therefore, based on the power-conserving interconnection, the controller
can be used to inject additional damping. As power-conserving feedback
structure we set
\[
\begin{array}{ccccccc}
\hat{u}_{c} & = & \hat{K}\rfloor\hat{y}, & \quad & \check{u}_{c} & = & \check{K}\rfloor\check{y},\\
\hat{u} & = & -\hat{K}^{*}\rfloor\hat{y}_{c}, & \quad & \check{u} & = & -\check{K}^{*}\rfloor\check{y}_{c},
\end{array}
\]
clearly satisfying (\ref{eq:power_conserving_interconnection}), together
with appropriate maps $\hat{K},\check{K}$ and their duals $\hat{K}^{*},\check{K}^{*}$,
respectively. It is worth stressing that the closed-loop system still
possesses a pH-structure with $\mathscr{H}_{cl}=\int_{\mathcal{B}}\mathcal{H}\mathrm{d}z^{1}+H_{c}$
as closed-loop Hamiltonian.

Next, we investigate Casimir functionals (structural invariants) of
the closed loop, which shall enable us to relate some of the controller
states to the plant in order to partially shape $\mathscr{H}_{cl}$.
Note that the controller states that are not related to the plant
can be used for the damping injection and thus for the purpose of
stabilisation.\begin{rem}\label{rem:stability}It should be noted
that in this contribution the focus is on a formal approach based
on differential-geometric methods. Thus no detailed stability investigations
will be carried out as this requires functional-analytic methods in
general. However, the relations $\mathscr{H}_{cl}>0$ and $\dot{\mathscr{H}}_{cl}\leq0$
serve as necessary conditions for a stability investigation in the
sense of Lyapunov. Worth stressing is the fact that for the nonlinear
PDE system under investigation, cf. Ex \ref{ex:nonlin_beam_diss},
the proof of stability is no trivial task as the verification of the
pre-compactness of the closed-loop trajectories is not straightforward.\end{rem}

\subsection{Determination of the Conditions for Structural Invariants}

Motivated by \cite{Siuka2011a}, we introduce the specific functionals
\[
\mathscr{C}^{\lambda}=x_{c}^{\lambda}+\int_{\mathcal{B}}\mathcal{C}^{\lambda}\mathrm{d}z^{1},\quad\mathcal{C}^{\lambda}\in C^{\infty}(\mathcal{J}^{2}(\mathcal{E})),
\]
with $\lambda=1,\ldots,\bar{n}_{c}\leq n_{c}$, which have to fulfil
$\dot{\mathscr{C}}^{\lambda}=0$ independently of $\mathcal{H}$ and
$H_{c}$ in order to serve as conserved quantities. That is, they
have to meet the conditions\begin{subequations}\label{eq:casimir_conditions}
\begin{align}
(J_{c}^{\lambda\beta_{c}}-R_{c}^{\lambda\beta_{c}}) & =0\\
\delta_{\alpha}\mathcal{C}^{\lambda}(\mathfrak{J}^{\alpha\beta}-\mathfrak{R}^{\alpha\beta}) & =0\label{eq:casimir_pdes}\\
(\hat{G}_{c,\hat{\xi}}^{\lambda}\hat{K}^{\hat{\xi}\hat{\eta}}\hat{B}_{\alpha\hat{\eta}}+\delta_{\alpha}^{\partial,1}\mathcal{C}^{\lambda})|_{\partial\mathcal{B}_{a}} & =0\\
(\check{G}_{c,\check{\xi}}^{\lambda}\check{K}^{\check{\xi}\check{\eta}}\check{B}_{\alpha\check{\eta}}+\delta_{\alpha}^{\partial,2}\mathcal{C}^{\lambda})|_{\partial\mathcal{B}_{a}} & =0\\
(\dot{x}^{\alpha}\delta_{\alpha}^{\partial,1}\mathcal{C}^{\lambda}+\dot{x}_{1}^{\alpha}\delta_{\alpha}^{\partial,2}\mathcal{C}^{\lambda})|_{\partial\mathcal{B}_{u}} & =0\label{eq:boundary_condition_casimir}
\end{align}
\end{subequations}to qualify as structural invariants. For the derivation
of (\ref{eq:casimir_conditions}) we refer to \cite[Eqs. (18) and (19)]{Rams2017a},
where it must be emphasised that only the non-differential operator
case is treated there. However, based on the chosen boundary-port
parameterisation, the computation can easily be adopted to the specific
differential-operator case considered here. Note that (\ref{eq:casimir_pdes})
can easily be satisfied by setting $\mathcal{C}^{\lambda}=d_{1}(\bar{\mathcal{C}}^{\lambda})$
with $\bar{\mathcal{C}}^{\lambda}\in C^{\infty}(\mathcal{J}^{1}(\mathcal{E}))$,
as the variational derivative annihilates total derivatives (see \cite[Theorem 4.7]{Olver1993}),
i.e. $\delta_{\alpha}\mathcal{C}^{\lambda}=0$ is met independently
of the concrete function $\bar{\mathcal{C}}^{\lambda}$.

\subsection{Energy-Casimir Controller for Example \ref{ex:nonlin_beam_diss}}

In \cite{Rams2017a}, an energy-based control law was developed for
the vertical deflection of a linear Euler-Bernoulli beam without dissipation.
The aim of this subsection is to design a Casimir-based controller
for the nonlinear beam structure with dissipation of Ex. \ref{ex:nonlin_beam_diss},
in order to stabilise the desired equilibrium (with arbitrary constants
$a,b\in\mathbb{R}$)
\begin{equation}
w^{1,d}=az^{1}+b,\quad w^{2,d}=0,\quad w_{1}^{1,d}=a.\label{eq:des_equilib}
\end{equation}
In the following, we design a nonlinear dynamical controller with
$n_{c}=6$, where it must be emphasised that three controller coordinates
are related to the plant to properly shape the Hamiltonian such that
(\ref{eq:des_equilib}) becomes a part of the minimum. The remaining
three controller states are used to inject additional damping in the
closed loop. If we consider the total derivatives $\mathcal{C}^{1}=-\frac{1}{L}d_{1}(z^{1}w^{1})$,
$\mathcal{C}^{2}=-\frac{1}{L}d_{1}(z^{1}w^{2})$ and $\mathcal{C}^{3}=-\frac{1}{L}d_{1}(z^{1}w_{1}^{1})$
as Casimir functions, we find that the conditions (\ref{eq:casimir_conditions})
yield the algebraic restrictions
\begin{align}
J_{c}-R_{c} & =\left[\begin{array}{cc}
0_{3\times3} & 0_{3\times3}\\
0_{3\times3} & A
\end{array}\right]\nonumber \\
\hat{G}_{c} & =\left[\begin{array}{cccccc}
1 & 0 & 0 & \hat{G}_{c,1}^{4} & \hat{G}_{c,1}^{5} & \hat{G}_{c,1}^{6}\\
0 & 1 & 0 & \hat{G}_{c,2}^{4} & \hat{G}_{c,2}^{5} & \hat{G}_{c,2}^{6}
\end{array}\right]\label{eq:controller_mappings}\\
\check{G}_{c} & =\left[\begin{array}{cccccc}
0 & 0 & 1 & \check{G}_{c,1}^{4} & \check{G}_{c,1}^{5} & \check{G}_{c,1}^{6}\end{array}\right]\nonumber 
\end{align}
with
\[
A=\left[\begin{array}{ccc}
-R_{c}^{44} & J_{c}^{45}-R_{c}^{45} & J_{c}^{46}-R_{c}^{46}\\
-J_{c}^{45}-R_{c}^{45} & -R_{c}^{55} & J_{c}^{56}-R_{c}^{56}\\
-J_{c}^{46}-R_{c}^{46} & -J_{c}^{56}-R_{c}^{56} & -R_{c}^{66}
\end{array}\right]
\]
for the controller mappings. Then the Casimir functions yield the
important relations
\[
\begin{array}{ccccc}
\mathscr{C}^{1} & = & x_{c}^{1}+\int_{0}^{L}\mathcal{C}^{1}\mathrm{d}z^{1} & = & x_{c}^{1}-w^{1}|_{L},\\
\mathscr{C}^{2} & = & x_{c}^{2}+\int_{0}^{L}\mathcal{C}^{2}\mathrm{d}z^{1} & = & x_{c}^{2}-w^{2}|_{L},\\
\mathscr{C}^{3} & = & x_{c}^{3}+\int_{0}^{L}\mathcal{C}^{3}\mathrm{d}z^{1} & = & x_{c}^{3}-w_{1}^{1}|_{L}.
\end{array}
\]
By choosing the initial conditions for the controller states appropriately,
the relations between the plant and the controller states result in
\[
\begin{array}{ccc}
x_{c}^{1}=w^{1}|_{L}, & x_{c}^{2}=w^{2}|_{L}, & x_{c}^{3}=w_{1}^{1}|_{L}\end{array}.
\]
To appropriately shape $\mathscr{H}_{cl}$, we set $H_{c}$ to
\begin{multline*}
H_{c}=\frac{c_{1}}{4}(x_{c}^{1}-w^{1,d}|_{L})^{4}+\frac{c_{2}}{4}(x_{c}^{2}-w^{2,d}|_{L})^{4}+\ldots\\
+\frac{c_{3}}{4}(x_{c}^{3}-w_{1}^{1,d}|_{L})^{4}+\frac{1}{2}M_{c,\mu\nu}x_{c}^{\mu}x_{c}^{\nu}
\end{multline*}
together with the constants $c_{1},c_{2},c_{3}>0$ and the positive
definite matrix $M_{c}$, $M_{c,\mu\nu}\in\mathbb{R}$ for $\mu,\nu=4,5,6$.
Consequently, $\mathscr{H}_{cl}$ evolves along solutions of the closed
loop according to 
\begin{multline*}
\dot{\mathscr{H}}_{cl}=-\int_{\mathcal{B}}(\alpha_{1}(\dot{w}_{11}^{1})^{2}+\alpha_{2}(\dot{w}_{1}^{2})^{2})\mathrm{d}z^{1}+\ldots\\
-x_{c}^{\mu}M_{c,\mu\nu}R_{c}^{\nu\rho}M_{c,\rho\vartheta}x_{c}^{\vartheta},
\end{multline*}
with $\rho,\vartheta=4,5,6$, from which it becomes apparent that
the controller injects additional damping. 

We want to stress that although the dynamic controller is able to
stabilise the desired rest position, there is no systematic approch
to determine the remaining degrees of freedom for the controller maps
properly. Therefore, in the next section, we propose a static controller
based on energy balancing that enables to simplify the controller
design compared to the energy-Casimir method.

\section{Energy-Balancing Control\label{sec:Energy-Balancing-Control}}

The intention of this section is to introduce a further control methodology
which is able to shape the plant Hamiltonian (at least partially)
and to add damping. In \cite{Macchelli2017}, a similar control framework
is proposed exploiting a pH system representation based on Stokes-Dirac
structures which is suitable for linear PDEs. Therefore, we propose
a similar control scheme for pH-systems with 2nd-order Hamiltonian
formulated in terms of the jet-bundle framework, which makes heavy
use of the geometric properties of the variational derivative as well
as both the boundary operators, and is not restricted to linear PDE
systems and controllers.

First, we demonstrate the basic principle of the energy-balancing
control (EBC) scheme by means of the non-differential operator case,
i.e. for systems of the form (\ref{eq:pH_non_differential_operator}).
Then, since the introduction of boundary ports is not straighforward
in the differential-operator case, the proposed control methodology
is studied on the example of the nonlinear beam structure with structural
damping of Ex. \ref{ex:nonlin_beam_diss}.

In the following, we consider static control laws that can be divided
into an energy-shaping part $\beta$ and a damping-injection part
$u'$ of the form 
\begin{equation}
\hat{u}^{\hat{\xi}}=\hat{\beta}^{\hat{\xi}}+\hat{u}'^{\hat{\xi}},\qquad\check{u}^{\check{\mu}}=\check{\beta}^{\check{\mu}}+\check{u}'^{\check{\mu}}.\label{eq:EBC_ansatz}
\end{equation}
The aim is to use the energy-shaping input $\beta$ to map the open-loop
equations (\ref{eq:pH_non_differential_operator}) into the target
system
\begin{equation}
\dot{x}^{\alpha}=\left(\mathcal{J}^{\alpha\beta}-\mathcal{R}^{\alpha\beta}\right)\delta_{\beta}\mathcal{H}_{d},\label{eq:target_System}
\end{equation}
with the ansatz $\mathcal{H}_{d}=\mathcal{H}+\mathcal{H}_{a}$, where
$\mathcal{H}_{a}$ is chosen such that $\mathcal{H}_{d}$ has a minimum
at the desired position of rest. Moreover, we set the (new) input
of the target system to
\begin{equation}
\begin{array}{cccccc}
\hat{B}_{\alpha\hat{\xi}}\hat{u}'^{\hat{\xi}} & = & \delta_{\alpha}^{\partial,1}\mathcal{H}_{d}|_{\partial\mathcal{B}_{a}}, & \check{B}_{\alpha\check{\mu}}\check{u}'^{\check{\mu}} & = & \delta_{\alpha}^{\partial,2}\mathcal{H}_{d}|_{\partial\mathcal{B}_{a}}\end{array}\label{eq:closed_loop_input}
\end{equation}
which shall be used to inject some additional damping.

In what follows, we explain in detail how the energy-shaping as well
as the damping-injection part have to be determined. It should be
noted that since we take no distributed input into account, we are
not able to modify the original dynamics of the system (\ref{eq:pH_non_differential_operator}).
Consequently, the matching equations are given as
\[
(\mathcal{J}^{\alpha\beta}-\mathcal{R}^{\alpha\beta})\delta_{\beta}\mathcal{H}=(\mathcal{J}^{\alpha\beta}-\mathcal{R}^{\alpha\beta})\delta_{\beta}\mathcal{H}_{d},
\]
and due to $\mathcal{H}_{d}=\mathcal{H}+\mathcal{H}_{a}$, they can
be reduced to
\begin{equation}
(\mathcal{J}^{\alpha\beta}-\mathcal{R}^{\alpha\beta})\delta_{\beta}\mathcal{H}_{a}=0,\label{eq:simpl_mc}
\end{equation}
yielding first conditions for $\mathcal{H}_{a}$. Furthermore, due
to the fact that there is no power flow through the unactuated boundary
$\partial\mathcal{B}_{u}$, $\mathcal{H}_{a}$ must satisfy the boundary
conditions $\delta_{\alpha}^{\partial,1}\mathcal{H}_{a}|_{\partial\mathcal{B}_{u}}=0$
and $\delta_{\alpha}^{\partial,2}\mathcal{H}_{a}|_{\partial\mathcal{B}_{u}}=0$,
cf. (\ref{eq:boundary_condition_casimir}). In the following, we make
use of the pleasant property that the variational derivative always
annihilates total derivatives, which allows us to derive a proper
control law in an elegant manner. In accordance with both boundary
categories, cf. (\ref{eq:output_param}), a suitable choice for the
additional Hamiltonian density is given as
\[
\mathcal{H}_{a}=\sum_{\hat{\nu}=1}^{\hat{m}}\hat{h}_{\hat{\nu}}+\sum_{\check{\mu}=1}^{\check{m}}\check{h}_{\check{\mu}},
\]
with $\hat{h}_{\hat{\nu}}=d_{1}(\hat{f}_{\hat{\nu}})$ and $\check{h}_{\check{\mu}}=d_{1}(\check{f}_{\check{\mu}})$
meeting the required boundary conditions. Thus, we are able to modify
the system's actuated boundary $\partial\mathcal{B}_{a}$ ensured
that (\ref{eq:simpl_mc}) is satisfied since $\delta_{\alpha}\mathcal{H}_{a}=0$
holds \textendash{} due to the special choice for $\hat{h}_{\hat{\nu}}$
and $\check{h}_{\check{\mu}}$ \textendash , which enables us to develop
nonlinear control laws as well. In fact, with regard to control purposes,
we are interested how the closed-loop functional $\mathscr{H}_{d}=\int_{\mathcal{B}}\mathcal{H}_{d}\mathrm{d}z^{1}$
evolves along solutions of the closed loop. A straightforward evaluation
yields
\begin{multline}
\dot{\mathscr{H}}_{d}=-\int_{\mathcal{B}}\delta_{\alpha}(\mathcal{H}_{d})\mathcal{R}^{\alpha\beta}\delta_{\beta}(\mathcal{H}_{d})\mathrm{d}z^{1}+\ldots\\
+(\dot{x}^{\alpha}\delta_{\alpha}^{\partial,1}\mathcal{H}_{d}+\dot{x}_{1}^{\alpha}\delta_{\alpha}^{\partial,2}\mathcal{H}_{d})|_{\partial\mathcal{B}_{a}},\label{eq:H_p_closed_loop}
\end{multline}
and by recalling that $\delta_{\alpha}\mathcal{H}_{a}=0$ is met,
we see that the dissipation effect in the domain remains unchanged.
If we consider the boundary expression of (\ref{eq:H_p_closed_loop}),
with (\ref{eq:EBC_ansatz}) and (\ref{eq:closed_loop_input}) we find
that the energy-shaping control laws result in
\begin{align}
\hat{B}_{\alpha\hat{\xi}}\hat{\beta}^{\hat{\xi}} & =-\delta_{\alpha}^{\partial,1}\mathcal{H}_{a}|_{\partial\mathcal{B}_{a}},\label{eq:EBC_control_law_a}\\
\check{B}_{\alpha\check{\xi}}\check{\beta}^{\check{\xi}} & =-\delta_{\alpha}^{\partial,2}\mathcal{H}_{a}|_{\partial\mathcal{B}_{a}}.\label{eq:EBC_control_law_b}
\end{align}
Now, since we mapped the original system into the target system (\ref{eq:target_System}),
we use its collocated boundary pairs $(\hat{u}'^{\hat{\eta}},\hat{y}_{\hat{\eta}})$
and $(\check{u}'^{\check{\mu}},\check{y}_{\check{\mu}})$ to inject
additional damping. If we choose the inputs for the damping-injection
part according to
\[
\hat{u}'^{\hat{\eta}}=-\hat{K}^{\hat{\eta}\hat{\kappa}}\hat{y}_{\hat{\kappa}},\quad\check{u}'^{\check{\mu}}=-\check{K}^{\check{\mu}\check{\rho}}\check{y}_{\check{\rho}},
\]
with $[\hat{K}^{\hat{\eta}\hat{\kappa}}]>0$ and $[\check{K}^{\check{\mu}\check{\rho}}]>0$,
then (\ref{eq:H_p_closed_loop}) follows to
\[
\dot{\mathscr{H}}_{d}=-\int_{\mathcal{B}}(\delta_{\alpha}\mathcal{H})\mathcal{R}^{\alpha\beta}(\delta_{\beta}\mathcal{H})\mathrm{d}z^{1}-\hat{y}_{\hat{\eta}}\hat{K}^{\hat{\eta}\hat{\kappa}}\hat{y}_{\hat{\kappa}}-\check{y}_{\check{\mu}}\check{K}^{\check{\mu}\check{\rho}}\check{y}_{\check{\rho}}
\]
ensuring that $\dot{\mathscr{H}}_{d}\leq0$ is met, i.e. that $\mathscr{H}_{d}$
is non-increasing along closed-loop solutions; cf. Remark \ref{rem:stability}
regarding the stability investigation.

\subsection{Energy-Balancing Controller for Example \ref{ex:nonlin_beam_diss}}

Next, we investigate the impact of linear differential operators to
the proposed controller design procedure. In particular, we develop
a control law in order to stabilise the desired equilibrium (\ref{eq:des_equilib})
for the nonlinear beam with structural damping of Ex. \ref{ex:nonlin_beam_diss}.
Due to the fact that the influence of $\mathfrak{R}_{A}$ and $\mathfrak{R}_{B}$
cannot be given in general, a parameterisation of the input of (\ref{eq:target_System})
according to (\ref{eq:closed_loop_input}) is not possible. However,
for the considered beam structure we set $\hat{u}'^{1}=\breve{Q}+\delta_{1}^{\partial,1}\mathcal{H}_{a}$,
$\hat{u}'^{2}=\breve{N}+\delta_{2}^{\partial,1}\mathcal{H}_{a}$ and
$\check{u}'^{1}=\breve{M}+\delta_{1}^{\partial,2}\mathcal{H}_{a}$;
therefore, the boundary maps result in $\hat{B}_{\alpha\hat{\xi}}=\delta_{\alpha\hat{\xi}}$
and $\check{B}_{\alpha\check{\xi}}=\delta_{\alpha\check{\xi}}$, together
with the Kronecker-Delta symbol meeting $\delta_{\alpha\beta}=1$
for $\alpha=\beta$ and $\delta_{\alpha\beta}=0$ for $\alpha\neq\beta$.

\begin{figure}
\input{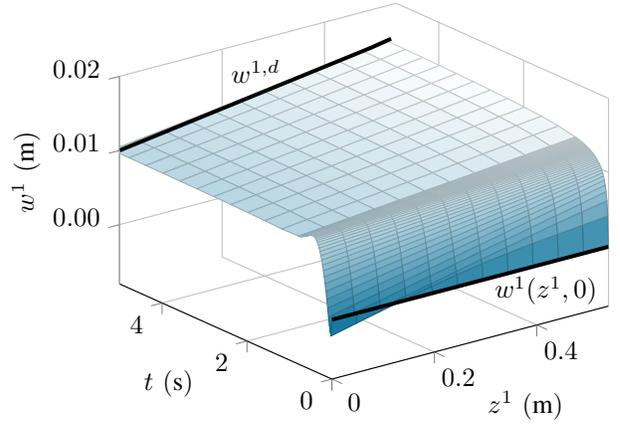}

\caption{\label{fig:EBC_nonlin_beam}Simulation results for $w^{1}$ for the
proposed EBC approach. Beam parameters: $L=0.54\,\mathrm{m}$, $EI=14.97\,\mathrm{Nm^{2}}$,
$EA=50\,\mathrm{N}$, $\rho A=2.1\,\mathrm{kgm^{-1}}$ and initial
rest position $w^{1}(z^{1},0)=0$. Controller parameters: $c_{1}=2\cdot10^{8}$,
$c_{2}=1000$, $c_{3}=8\cdot10^{4}$, $k_{1}=2200$ and $k_{2}=k_{3}=1$.}
\end{figure}

Following the intention of the proposed control scheme, we choose
the additional Hamiltonian according to $\mathcal{H}_{a}=\hat{h}_{1}+\hat{h}_{2}+\check{h}_{1}$,
with
\[
\begin{array}{cccccc}
\hat{h}_{1} & = & d_{1}(\hat{f}_{1}),\quad & \hat{f}_{1} & = & \frac{z^{1}}{4L}c_{1}(w^{1}-w^{1,d})^{4},\\
\hat{h}_{2} & = & d_{1}(\hat{f}_{2}),\quad & \hat{f}_{2} & = & \frac{z^{1}}{4L}c_{2}(w^{2}-w^{2,d})^{4},\\
\check{h}_{1} & = & d_{1}(\check{f}_{1}),\quad & \check{f}_{1} & = & \frac{z^{1}}{4L}c_{3}(w_{1}^{1}-w_{1}^{1,d})^{4},
\end{array}
\]
in order to shape the closed-loop Hamiltonian at the actuated boundary
$\partial\mathcal{B}_{a}$ and to fulfil the boundary conditions $\delta_{\alpha}^{\partial,1}\mathcal{H}_{a}|_{\partial\mathcal{B}_{u}}=0$
and $\delta_{\alpha}^{\partial,2}\mathcal{H}_{a}|_{\partial\mathcal{B}_{u}}=0$.
An evaluation of (\ref{eq:EBC_control_law_a}) and (\ref{eq:EBC_control_law_b})
yields the EBC laws
\begin{align*}
\hat{\beta}^{1} & =-c_{1}(w^{1}|_{L}-w^{1,d}|_{L})^{3},\\
\hat{\beta}^{2} & =-c_{2}(w^{2}|_{L}-w^{2,d}|_{L})^{3},\\
\check{\beta}^{1} & =-c_{3}(w_{1}^{1}|_{L}-w_{1}^{1,d}|_{L})^{3},
\end{align*}
mapping the beam into the target system $\dot{x}^{\alpha}=(\mathcal{J}^{\alpha\beta}-\mathfrak{R}^{\alpha\beta})\delta_{\beta}\mathcal{H}_{d}$
with the Hamiltonian density $\mathcal{H}_{d}=\mathcal{H}+\mathcal{H}_{a}$.
If we set the (damping-injection) input for this system to
\[
\begin{array}{cccccc}
\hat{u}'^{1} & = & -k_{1}(\dot{w}^{1}|_{L}), & \hat{u}'^{2} & = & -k_{2}(\dot{w}^{2}|_{L}),\\
\check{u}'^{1} & = & -k_{3}(\dot{w}_{1}^{1}|_{L}),
\end{array}
\]
the power-balance relation results in
\begin{multline*}
\dot{\mathscr{H}}_{d}=-\int_{\mathcal{B}}(\alpha_{1}(\dot{w}_{11}^{1})^{2}+\alpha_{2}(\dot{w}_{1}^{2})^{2})\mathrm{d}z^{1}+\ldots\\
-k_{1}(\dot{w}^{1}|_{L})^{2}-k_{2}(\dot{w}^{2}|_{L})^{2}-k_{3}(\dot{w}_{1}^{1}|_{L})^{2},
\end{multline*}
obviously satisfying $\dot{\mathscr{H}}_{d}\leq0$. It should be noted
that the EBC scheme suggests a nonlinear PD control law, where the
proportional part may be interpreted as a nonlinear spring. From Figure
\ref{fig:EBC_nonlin_beam}, it becomes clearly evident that the proposed
controller stabilises the desired equilibrium (\ref{eq:des_equilib})
with $a=0.01$ and $b=0.01$.

\section{Conclusion\label{sec:Conclusion}}

In this contribution, the pH-framework based on the jet-bundle scenario
has been used for the modeling and energy-based control for a nonlinear
Euler-Bernoulli beam including structural damping. In particular,
the control by interconnection method has been adapted for the system
class under investigation and for the first time, EBC techniques that
have been already considered in the Stokes-Dirac scenario (linear
setting), have been analysed for nonlinear PDEs based on the geometric
framework exploiting the underlying jet-bundle structure. Future research
directions will include the generalisation of EBC techniques to domain
inputs in this setting and extensions towards higher-dimensional spatial
domains.

\section*{Appendix}

In the appendix, two examples are studied that serve as a motivation
for the presented framework in the non-differential operator case.
In particular, a linear and a nonlinear Euler-Bernoulli beam are considered
where in contrast to Ex. \ref{ex:nonlin_beam_diss} no differential
operators appear.\begin{exmp}[Linear beam structure with viscous damping]\label{ex:lin_beam}First,
we investigate a beam with linearised geometric and linear constitutive
relations, combined with external viscous damping effects. The corresponding
equations of motion in classical notation read as
\begin{equation}
\begin{array}{ccc}
\rho A\ddot{w}^{1} & = & -EI\frac{\partial^{4}w^{1}}{\partial\left(z^{1}\right)^{4}}-\alpha_{1}\dot{w}^{1},\\
\rho A\ddot{w}^{2} & = & EA\frac{\partial^{2}w^{2}}{\partial\left(z^{1}\right)^{2}}-\alpha_{2}\dot{w}^{2},
\end{array}\label{eq:lin_beam_visc_damp}
\end{equation}
with $E,I,\rho,A>0$ as material parameters, see \cite{Meirovitch1997}
for more details, and small positive constants $\alpha_{1},\alpha_{2}$
for the damping parameters. Because we are interested in a system
description according to (\ref{eq:pH_non_differential_operator}),
we first introduce the generalised momenta $p_{1}=\rho A\dot{w}^{1}$
and $p_{2}=\rho A\dot{w}^{2}$. By means of these new coordinates,
the total energy density of the beam can be expressed by
\[
\mathcal{H}=\frac{1}{2\rho A}((p_{1})^{2}+(p_{2})^{2})+\frac{1}{2}EI(w_{11}^{1})^{2}+\frac{1}{2}EA(w_{1}^{2})^{2}.
\]
Consequently, we find that (\ref{eq:lin_beam_visc_damp}) can be formulated
as
\begin{equation}
\left[\begin{array}{c}
\begin{array}{c}
\dot{w}^{1}\\
\dot{w}^{2}
\end{array}\\
\dot{p}_{1}\\
\dot{p}_{2}
\end{array}\right]=(\mathcal{J}-\mathcal{R})\left[\begin{array}{c}
\delta_{w^{1}}\mathcal{H}\\
\delta_{w^{2}}\mathcal{H}\\
\delta_{p_{1}}\mathcal{H}\\
\delta_{p_{2}}\mathcal{H}
\end{array}\right],\label{eq:pH_lin_beam}
\end{equation}
where the tensors $\mathcal{J}$ and $\mathcal{R}$ correspond to
\[
\mathcal{J}=\left[\begin{array}{cccc}
0 & 0 & 1 & 0\\
0 & 0 & 0 & 1\\
-1 & 0 & 0 & 0\\
0 & -1 & 0 & 0
\end{array}\right],\quad\mathcal{R}=\left[\begin{array}{cccc}
0 & 0 & 0 & 0\\
0 & 0 & 0 & 0\\
0 & 0 & \alpha_{1} & 0\\
0 & 0 & 0 & \alpha_{2}
\end{array}\right].
\]
 A straightforward evaluation of (\ref{eq:h_p_non_diff_op}) yields
\begin{multline*}
\dot{\mathscr{H}}=-\int_{\mathcal{B}}(\alpha_{1}(\dot{w}^{1})^{2}+\alpha_{2}(\dot{w}^{2})^{2})\mathrm{d}z^{1}+\ldots\\
+(\dot{w}^{1}Q+\dot{w}^{2}N+\dot{w}_{1}^{1}M)|_{\partial\mathcal{B}_{a}},
\end{multline*}
with the shear force $Q=-EIw_{111}^{1}$, the normal force $N=EAw_{1}^{2}$
and the bending moment $M=EIw_{11}^{1}$, where we made extensive
use of (\ref{eq:boundary_operators}). Because the part of $\mathcal{H}$
which is associated to the transveral vibrations is of 2nd order,
we obtain the two boundary ports $(\dot{w}^{1},Q)$ and $(\dot{w}_{1}^{1},M)$.
On the other hand, the part of $\mathcal{H}$ which is related to
the longitudinal vibrations of the beam is of first order, and thus
only one boundary port $(\dot{w}^{2},N)$ appears. It is worth stressing
that these three ports are real power ports because the Hamiltonian
corresponds to the total energy of the beam. Note that if $Q$, $N$
and $M$ serve as controlled variables this fixes the corresponding
collocated outputs, and the construction of the maps $\hat{B}_{\alpha\hat{\xi}}$
and $\check{B}_{\alpha\check{\xi}}$ of (\ref{eq:output_param}) becomes
straightforward in that case.\end{exmp}

\begin{exmp}[Nonlinear 2nd-order beam structure]\label{ex:nonlin_beam}Next,
we consider nonlinear geometric relations within the Euler-Bernoulli
beam theory. In particular, we consider nonlinear strain-deflection
relations in the sense of von Kármán, see \cite{Kugi2001} for example.
Using this strain-deflection relations leads to the potential energy
density
\[
\mathcal{P}=\frac{1}{2}EA((w_{1}^{2})^{2}+\frac{1}{4}(w_{1}^{1})^{4}+w_{1}^{2}(w_{1}^{1})^{2})+\frac{1}{2}EI(w_{11}^{1})^{2},
\]
while for the kinetic energy density we have $\mathcal{K}=\frac{1}{2}\rho A\left((w^{1})^{2}+(w^{2})^{2}\right)$.
Then, applying Hamilton's principle (see \cite{Meirovitch1997}) on
the Lagrangian density $\mathcal{L}=\mathcal{K}-\mathcal{P}$ yields
the following equations of motion 
\begin{align}
\ddot{w}^{1} & =\frac{E}{\rho}(\frac{3}{2}w_{11}^{1}(w_{1}^{1})^{2}+w_{11}^{2}w_{1}^{1}+w_{1}^{2}w_{11}^{1})-\frac{EI}{\rho A}w_{1111}^{1},\nonumber \\
\ddot{w}^{2} & =\frac{E}{\rho}w_{11}^{2}+\frac{E}{\rho}w_{1}^{1}w_{11}^{1}.\label{eq:nonlin_beam_eom}
\end{align}
Again, by using the generalised momenta $p_{1}=\rho A\dot{w}^{1}$
and $p_{2}=\rho A\dot{w}^{2}$, as well as the total energy density
$\mathcal{K}+\mathcal{P}$ as Hamiltonian density, we find that the
governing equations can also be formulated as
\begin{equation}
\left[\begin{array}{c}
\begin{array}{c}
\dot{w}^{1}\\
\dot{w}^{2}
\end{array}\\
\dot{p}_{1}\\
\dot{p}_{2}
\end{array}\right]=\left[\begin{array}{cccc}
0 & 0 & 1 & 0\\
0 & 0 & 0 & 1\\
-1 & 0 & 0 & 0\\
0 & -1 & 0 & 0
\end{array}\right]\left[\begin{array}{c}
\delta_{w^{1}}\mathcal{H}\\
\delta_{w^{2}}\mathcal{H}\\
\delta_{p_{1}}\mathcal{H}\\
\delta_{p_{2}}\mathcal{H}
\end{array}\right].\label{eq:pH_nonlin_Beam}
\end{equation}
Note that (\ref{eq:pH_lin_beam}) and (\ref{eq:pH_nonlin_Beam}) are
structurally equivalent apart from the dissipation effects in (\ref{eq:pH_lin_beam});
however, the nonlinearity of (\ref{eq:pH_nonlin_Beam}) is hidden
in $\mathcal{H}$. Moreover, since no dissipation effects have been
taken into account, the formal change of the Hamiltonian functional
reads as 
\begin{equation}
\dot{\mathscr{H}}=(\dot{w}^{1}\tilde{Q}+\dot{w}^{2}\tilde{N}+\dot{w}_{1}^{1}\tilde{M})|_{\partial\mathcal{B}_{a}},\label{eq:h_p_nonlin_beam}
\end{equation}
with the modified shear force, normal force and bending moment that
follow again by using (\ref{eq:boundary_operators}) according to
\begin{align}
\tilde{Q} & =\frac{1}{2}EA(w_{1}^{1})^{3}+EAw_{1}^{2}w_{1}^{1}-EIw_{111}^{1},\nonumber \\
\tilde{N} & =EAw_{1}^{2}+\frac{1}{2}EA(w_{1}^{1})^{2},\label{eq:nonlin_beam_internal_forces}\\
\tilde{M} & =EIw_{11}^{1},\nonumber 
\end{align}
respectively. \end{exmp}

\bibliographystyle{ieeetr}
\bibliography{my_bib}

\begin{thebibliography}{10}

\bibitem{Schaft2000}
A.~J. van~der Schaft, {\em {L2-Gain and Passivity Techniques in Nonlinear
  Control}}.
\newblock {Springer}, 2000.

\bibitem{Ortega2001}
R.~Ortega, A.~J. van~der Schaft, I.~Mareels, and B.~Maschke, ``Putting energy
  back in control,'' {\em IEEE Control Syst. Mag.}, vol.~21, no.~2, pp.~18--33,
  2001.

\bibitem{Schaft2002}
A.~J. van~der Schaft and B.~Maschke, ``{Hamiltonian formulations of distributed
  parameter systems with boundary energy flow},'' {\em {Journal of Geometry and
  Physics}}, vol.~42, no.~1-2, pp.~166--194, 2002.

\bibitem{Gorrec2005}
Y.~L. Gorrec, H.~J. Zwart, and B.~Maschke, ``Dirac structures and boundary
  control systems associated with skew-symmetric differential operators,'' {\em
  {SIAM J. Control Optim.}}, vol.~44, no.~5, pp.~1864--1892, 2005.

\bibitem{Jacob2012}
B.~Jacob and H.~J. Zwart, {\em {Linear Port-Hamiltonian Systems on
  Infinite-dimensional Spaces}}.
\newblock {Birkh{\"a}user}, 2012.

\bibitem{Schoeberl2014a}
M.~Sch{\"o}berl and A.~Siuka, ``{Jet bundle formulation of infinite-dimensional
  port-Hamiltonian systems using differential operators},'' {\em {Automatica}},
  vol.~50, no.~2, pp.~607--613, 2014.

\bibitem{Schoeberl2015e}
M.~Sch{\"o}berl and K.~Schlacher, ``{Lagrangian and Port-Hamiltonian
  formulation for Distributed-parameter systems},'' in {\em {Proceedings of the
  8th Vienna International Conference on Mathematical Modelling}}, vol.~48,
  issue 1 of {\em {IFAC-PapersOnLine}}, pp.~610--615, 2015.

\bibitem{Curtain1995}
R.~F. Curtain and H.~J. Zwart, {\em {An Introduction to Infinite-Dimensional
  Linear System Theory}}.
\newblock {Springer}, 1995.

\bibitem{Macchelli2004}
A.~Macchelli and C.~Melchiorri, ``{Modeling and control of the Timoshenko beam.
  The distributed port Hamiltonian approach},'' {\em {SIAM J. Control Optim.}},
  vol.~43, no.~2, pp.~743--767, 2004.

\bibitem{Siuka2011a}
A.~Siuka, M.~Sch{\"o}berl, and K.~Schlacher, ``{Port-Hamiltonian modelling and
  energy-based control of the Timoshenko beam - An approach based on structural
  invariants},'' {\em {Acta Mechanica}}, vol.~222, pp.~69--89, 11 2011.

\bibitem{Rams2017a}
H.~Rams and M.~Sch{\"o}berl, ``{On Structural Invariants in the Energy Based
  Control of Port-Hamiltonian Systems with Second-Order Hamiltonian},'' in {\em
  {Proceedings of the American Control Conference (ACC)}}, pp.~1139--1144,
  2017.

\bibitem{Macchelli2017}
A.~Macchelli, Y.~L. Gorrec, H.~Ramirez, and H.~Zwart, ``{On the Synthesis of
  Boundary Control Laws for Distributed Port-Hamiltonian Systems},'' {\em {IEEE
  Trans. Autom. Control}}, vol.~62, no.~4, pp.~1700--1713, 2017.

\bibitem{Giachetta1997}
G.~Giachetta, L.~Mangiarotti, and G.~Sardanashvily, {\em {New Lagrangian and
  Hamiltonian Methods in Field Theory}}.
\newblock {World Scientific}, 1997.

\bibitem{Schoeberl2008a}
M.~Sch{\"o}berl, H.~Ennsbrunner, and K.~Schlacher, ``{Modelling of
  piezoelectric structures - a hamilton approach},'' {\em {Mathematical and
  Computer Modelling of Dynamical Systems}}, vol.~14, no.~3, pp.~179--193,
  2008.

\bibitem{Liu1998}
K.~Liu and Z.~Liu, ``{Exponential decay of energy of the Euler-Bernoulli beam
  with locally distributed Kelvin-Voigt damping},'' {\em {SIAM J. Control
  Optim.}}, vol.~36, no.~3, pp.~1086--1098, 1998.

\bibitem{Meirovitch1997}
L.~Meirovitch, {\em {Analytical Methods in Vibrations}}.
\newblock {Macmillan Publishing Co., Inc.}, 1967.

\bibitem{Kugi2001}
A.~Kugi, {\em {Non-linear Control Based on Physical Models: Electrical,
  mechanical and hydraulic systems}}, vol.~260 of {\em {Lecture Notes in
  Control and Information Sciences}}.
\newblock Springer, 2001.

\bibitem{Schoeberl2013a}
M.~Sch{\"o}berl and A.~Siuka, ``{On Casimir Functionals for
  infinite-dimensional Port-Hamiltonian Control Systems},'' {\em IEEE
  Transactions on Automatic Control}, vol.~58, no.~7, pp.~1823--1828, 2013.

\bibitem{Olver1993}
P.~J. Olver, {\em {Applications of Lie Groups to Differential Equations}}.
\newblock {Springer}, 2nd~ed., 1993.

\end{thebibliography}

\end{document}